\documentclass[twoside]{article}
\title{A Note on Gorenstein
Flat Dimension}
\date{}

\usepackage{amsfonts}
\usepackage{amsmath}
\usepackage{amssymb}
\usepackage{latexsym}

\newtheorem{thm}{\bf Theorem}[section]
\newtheorem{cor}[thm]{\bf Corollary}
\newtheorem{lem}[thm]{\bf Lemma}
\newtheorem{prop}[thm]{\bf Proposition}
\newtheorem{defn}[thm]{\bf Definition}

\catcode`\ç=13
\defç{\c{c}}
\catcode`\é=13
\defé{\'e}
\catcode`\à=13
\defà{\`a}
\catcode`\è=13
\defè{\`e}
\catcode`\â=13
\defâ{\^a}
\catcode`\ù=13
\defù{\`u}
\catcode`\ê=13
\defê{\^e}
\catcode`\î=13
\defî{\^\i}
\catcode`\ô=13
\defô{\^o}
\newcommand{\field}[1]{\mathbb{#1}}

\newcommand{\Q }{\field{Q}}
\newcommand{\Z }{\field{Z}}
\newcommand{\N }{\field{N}}

\def\proof{{\parindent0pt {\textit{Proof}.\ }}}

\def\fd{{\rm fd}}
\def\id{{\rm id}}

\def\Gpd{{\rm Gpd}}
\def\Gfd{{\rm Gfd}}
\def\Gid{{\rm Gid}}

\def\cfd{{\rm  cfd}}

\def\Im{{\rm Im}}

\def\Ext{{\rm Ext}}
\def\Tor{{\rm Tor}}
\def\Hom{{\rm Hom}}

\def\sup{{\rm sup}}

\begin{document}
\thispagestyle{empty}

\maketitle \vspace*{-2cm}

\begin{center}
\begin{large} \textbf{Driss Bennis}\end{large}\\[0.1cm]
\textsl{Department of Mathematics, Faculty of Science and Technology of
Fez,\\ Box 2202, University S. M. Ben Abdellah Fez,
Morocco\\[0.1cm]
E-mail: driss\_bennis@hotmail.com}
\end{center}
\noindent{\large\bf Abstract.} Unlike the Gorenstein projective and
injective dimensions, the majority of results on the Gorenstein flat
dimension have been established only over Noetherian (or coherent)
rings. Naturally, one would like to generalize these results to any
associative ring. In this direction, we show that the Gorenstein
flat dimension is a refinement of the classical flat dimension over
any ring; and we investigate the relations between the Gorenstein
projective dimension and the Gorenstein flat dimension.
\bigskip

 \medskip


\begin{section}{Introduction}
Throughout this paper, $R$ denotes an associative ring with
identity and all modules are unitary. For an  $R$-module $M$, we
use  $\id_R(M)$ and $\fd_R(M)$ to denote, respectively, the
classical injective and flat dimensions of $M$. We use $M^*$ to
denote the character module  $\Hom_{\Z}(M, \Q/\Z)$ of $M$.\bigskip

For every module over an associative ring, Enochs, Jenda, and
Torrecillas \cite{GoPlat} defined the Gorenstein flat dimension
(see Definition \ref{def-Gflat}) to complete the analogy between
 classical homological dimension theory and  Gorenstein
homological dimension theory. But, they mainly studied it when the
base ring is Gorenstein (see also \cite{Gf-cover}). Their
characterization of the Gorenstein flat dimension over Gorenstein
rings was generalized, by Chen and Ding \cite{FCring}, to  $n$-FC
rings (i.e., coherent rings with self-FP-injective dimension at
most $n$). Also, Christensen gave a characterization of the
Gorenstein flat dimension over local Cohen-Macaulay  rings with a
dualizing module \cite[Theorem 5.2.14]{LW} (see also \cite{Fox}).
Namely, in Christensen's book, it is shown that there are good
results for the Gorenstein flat dimension  over Noetherian rings,
which are very often local Cohen-Macaulay with a dualizing module.
In \cite{HH}, Holm generalized these results to coherent rings
(see also \cite{CFH}). In this paper, we give some results on the
Gorenstein flat dimension that hold over a larger class of
rings.\bigskip

In  Section 2, we  show that  the Gorenstein flat dimension is a
refinement of the usual flat dimension over any associative ring,
that is Theorem \ref{thm-Gf-f}. Then, in Theorem \ref{thm-Gf-Cf},
we establish a relation between the Gorenstein flat dimension and
the copure flat dimension (see Definition
\ref{def-copure}).\bigskip

In  Section  3, we investigate the relations between the Gorenstein
flat dimension and the Gorenstein projective dimension. Namely, we
establish another situation where the question ``Is every Gorenstein
projective module Gorenstein flat'' has an affirmative answer (see
Proposition \ref{prop-1-Gf-Gp}). In Theorem \ref{thm-Gp-Gf}, we show
that, if $M$ is an infinitely presented left $R$-module, then $M$ is
Gorenstein flat if  and only if  it is Gorenstein projective.

\end{section}


\begin{section}{Gorenstein flat, flat, and copure flat dimensions}
In this section, we investigate the relations  between the
Gorenstein flat dimension, the flat dimension, and the copure flat
dimension.\\
Recall the definitions of Gorenstein flat modules and Gorenstein
flat dimension.

\begin{defn}[\cite{LW}]\label{def-Gflat}
A complete  flat resolution is an exact sequence of flat  left
$R$-modules,
$$\mathbf{F}=\ \cdots\rightarrow
F_1\rightarrow F_0 \rightarrow F^0 \rightarrow F^1
\rightarrow\cdots,$$ such that
$I\otimes_R \mathbf{F}$ is exact for every injective right $R$-module $I$.\\
\indent A left $R$-module $M$ is called Gorenstein flat (G-flat
for short), if there exists a complete flat
resolution $\mathbf{F}$ with $M \cong \Im(F_0 \rightarrow F^0)$.\\
\indent For a positive integer $n$, we say that a left $R$-module
$M$ has Gorenstein flat dimension at most $n$, and we write
$\Gfd_R(M)\leq n$, if if there exists an exact sequence of
$R$-modules $  0\rightarrow G_n\rightarrow\cdots \rightarrow G_0
\rightarrow M \rightarrow 0,$ where each $G_i$ is Gorenstein flat.
\end{defn}

It is well-known that the Gorenstein homological dimensions are,
over Noetherian rings, refinements of the classical homological
dimensions. In  \cite{Rel-hom} and also in \cite{HH}, it is proved
that the result remains true for the Gorenstein projective and
injective dimensions over associative rings. However, in the note
after Proposition 3.6 of \cite{BM}, it is remarked that the
Gorenstein flat dimension is a refinement of the classical flat
dimension over coherent rings (see also \cite[Proposition 3.7 and
Corollary 3.8]{BM}). Next we show this holds over associative rings.

\begin{thm}\label{thm-Gf-f}
For a left $R$-module $M$, $\Gfd_R(M)\leq \fd_R(M)$ with equality
when $\fd_R(M)$ is finite.
\end{thm}
\proof The inequality is well-known and follows from the trivial
fact that every  flat module is Gorenstein flat (see for instance
\cite[Theorem 3.19]{HH}).\\
\indent Now assume that $\fd_R(M)<\infty$. Then, by   \cite[Lemma
3.51 and Theorem 3.52]{Rot}, $\id_R(M^*)=\fd_R(M)<\infty$ and so
$\Gid_R(M^*)=\id_R(M^*)$ (from the injective counterpart of
\cite[Proposition 2.27]{HH}). Therefore, combine the last equalities
and the inequality $\Gid_R(M^*)\leq \Gfd_R(M)$ of \cite[Theorem
3.11]{HH}, we get:
$$\fd_R(M)=\id_R(M^*)=\Gid_R(M^*)\leq \Gfd_R(M).$$
This gives the desired equality. \hfill$\square$\bigskip

Another dimension which is closely related with the Gorenstein
flat dimension is the copure flat dimension, which is defined as
follows:

\begin{defn}[\cite{copure}]\label{def-copure}
 The copure flat dimension of  a left $R$-module $M$, $\cfd_R(M)$, is defined
to be the largest positive integer $n$ such that
$\Tor_n^R(E,M)\not =0$ for some injective right $R$-module $E$.
\end{defn}

The copure flat dimension is used in the functorial description of
the Gorenstein flat dimension as follows \cite[Theorem 3.14]{HH}:
If $R$ is right coherent and $M$ is a left $R$-module with finite
Gorenstein flat dimension, then $\Gfd_R(M)=\cfd_R(M)$. This is a
generalization of \cite[Corollary 11]{FCring} and \cite[Theorem
5.2.14]{LW}. Over any ring we give the following:

\begin{thm}\label{thm-Gf-Cf} For any left $R$-module $M$, we have
inequality: $$\cfd_R(M)\leq \Gfd_R(M).$$ If $\fd_R(M)$ is finite,
then we have equalities:
$$\cfd_R(M)= \Gfd_R(M)=\fd_R(M).$$
\end{thm}

To prove this we need the following lemma.

\begin{lem}\label{lem-copure-flat} Let $M$ be a left $R$-module. If
$\fd_R(M)<\infty$, then $\cfd_R(M)= \fd_R(M).$
\end{lem}
\proof Obviously $\cfd_R(M)\leq \fd_R(M).$\\
\indent For the converse, suppose that $\fd_R(M)=n$ for some
positive integer $n$. Then there is  a right $R$-module   $N$ such
that $\Tor_n^R(N,M) \not = 0$, and it follows, by the long exact
sequence of Tor, that also $\Tor_n^R(E(N),M) \not = 0$, where $E(N)$
is the injective envelope of $N$. This implies that $\cfd_R(M)\geq
n=\fd_R(M)$, as desired. \hfill$\square$ \bigskip

\noindent\textit{Proof of Theorem \ref{thm-Gf-Cf}.} The inequality
$\cfd_R(M)\leq \Gfd_R(M) $ follows by dimension shifting argument
and using the fact that every Gorenstein flat left $R$-module  $G$
satisfies $\Tor_i^R(I,G)=0$ for all $i\geq 1$ and all injective
right $R$-modules
$I$.\\
\indent The equalities $\cfd_R(M)= \Gfd_R(M)=\fd_R(M) $ when
$\fd_R(M)$ is finite follow immediately by the inequalities
$\cfd_R(M) \leq \Gfd_R(M)\leq\fd_R(M) $  and Lemma
\ref{lem-copure-flat}. \hfill$\square$\bigskip

Note that the proof of the equalities of Theorem \ref{thm-Gf-Cf}
can be used as a proof of the equality of Theorem \ref{thm-Gf-f},
and so the two results can be written in one theorem. Here we
separate them because each one has a different aim.

\end{section}


\begin{section}{Gorenstein flat and Gorenstein projective dimensions}
In this section, we investigate the relations  between the
Gorenstein flat dimension and the  Gorenstein projective dimension.\\
First, recall the definitions of   Gorenstein projective modules
and   Gorenstein projective dimension.

\begin{defn}[\cite{LW}]\label{def-Gproj}
A complete  projective resolution is an exact sequence of
projective left $R$-modules,
$$\mathbf{P}=\ \cdots\rightarrow P_1\rightarrow P_0
\rightarrow P^0 \rightarrow P^1 \rightarrow\cdots,$$ such that
$\Hom_R (\mathbf{P}, Q) $ is exact for every projective left $R$-module $Q$.\\
\indent A left $R$-module $M$ is called Gorenstein projective
(G-projective for short), if there exists a complete projective
resolution $\mathbf{P}$ with $M \cong \Im(P_0 \rightarrow P^0)$.\\
\indent For a positive integer $n$, we say that a left $R$-module
$M$ has Gorenstein  projective dimension at most $n$, and we write
$\Gpd_R(M)\leq n$, if there exists an exact sequence of
$R$-modules $  0\rightarrow G_n\rightarrow\cdots \rightarrow G_0
\rightarrow M \rightarrow 0,$ where each $G_i$ is Gorenstein
projective.
\end{defn}

Motivated by the results in the classical case, there are two
principal questions concerning the relations between the
Gorenstein projective dimension and the Gorenstein flat dimension
of modules. In fact,  between the Gorenstein projective modules
and the Gorenstein flat modules:\bigskip

\noindent\textbf{Question A.} Is every Gorenstein projective
module Gorenstein flat?\\
\noindent\textbf{Question B.} When is a Gorenstein flat module
Gorenstein projective?\bigskip

Several attempts have been made to obtain, as the classical case,
an affirmative answer to Question A. In \cite[Proposition
3.4]{HH}, Holm proved, over right coherent ring with finite
finitistic projective dimension, that every Gorenstein projective
left  module is Gorenstein flat. In fact, this result holds, by
the proof of \cite[Proposition 3.4]{HH}, over right coherent rings
such that every flat left module has finite projective dimension,
which is remarked in \cite[Proposition 3.7]{CFH}. Here we give the
following result:

\begin{prop} \label{prop-1-Gf-Gp} If every injective right $R$-module has
finite flat dimension, then $\Gfd_R(M)\leq \Gpd_R(M)$ for every
left $R$-module $M$.
\end{prop}
\proof To prove the inequality, it is sufficient to prove that
every Gorenstein projective left module is Gorenstein flat, which
is equivalent to  proving that every complete projective
resolution
is complete flat.\\
Consider then a complete projective resolution $\mathbf{P}$ and an
injective right $R$-module $I$. By hypothesis, $\fd_R(I)=n$ for
some positive integer $n$. Then, we have an exact sequence
$$ 0 \rightarrow F_n\rightarrow \cdots
\rightarrow F_0\rightarrow I \rightarrow 0,$$ where each $F_i$ is
a flat right $R$-module. Let $I_{i}=\Im (F_i \rightarrow F_{i-1})$
for
$1\leq i \leq n-1$ and $I_0=I$.\\
Consider the short exact sequence $$ 0 \rightarrow F_n\rightarrow
F_{n-1}\rightarrow I_{n-1} \rightarrow 0.$$
Then,
 $$ 0 \rightarrow
F_n\otimes_R \mathbf{P}\rightarrow F_{n-1}\otimes_R
\mathbf{P}\rightarrow I_{n-1} \otimes_R \mathbf{P}\rightarrow 0$$
is a short exact sequence of complexes. Since $F_n$ and $F_{n-1}$
are flat,  $F_n\otimes_R \mathbf{P}$ and  $F_{n-1}\otimes_R
\mathbf{P}$
are exact, so is $I_{n-1} \otimes_R \mathbf{P}$.\\
Now consider the short exact sequence $$ 0 \rightarrow
I_{n-1}\rightarrow F_{n-2}\rightarrow I_{n-2} \rightarrow 0.$$
Then,
 $$ 0 \rightarrow
I_{n-1}\otimes_R \mathbf{P}\rightarrow F_{n-2}\otimes_R
\mathbf{P}\rightarrow I_{n-2} \otimes_R \mathbf{P}\rightarrow 0$$
is a short exact sequence of complexes. Since    $F_{n-2}$ is flat
and by the argument above, $I_{n-1}\otimes_R \mathbf{P}$ and
$F_{n-2}\otimes_R \mathbf{P}$
are exact, so is $I_{n-2} \otimes_R \mathbf{P}$.\\
\indent The argument above can be applied successively until we
conclude that the sequence $I \otimes_R \mathbf{P}$ is exact. This
implies the desired result. \hfill$\square$\bigskip

Note that the rings that satisfy the condition of Proposition
\ref{prop-1-Gf-Gp} were investigated in \cite{CD}.\bigskip

Now, we investigate Question B. In \cite{BM}, Question B was
investigated for a particular case of Gorenstein projective and
flat modules, such that well-known results in the classical case
were extended (please see \cite[Propositions 3.9 and 3.12 and
Corollary 3.10]{BM}). Here, we give a situation where a Gorenstein
flat module is Gorenstein projective.\bigskip

It is well-known that a finitely presented module is flat if and
only if it is projective. In Gorenstein homological dimension
theory, there is an analogous (in fact a generalization) of this
result over Noetherian rings \cite[Theorem 5.1.11]{LW}. Namely, it
is proved, over Noetherian rings, that a finitely generated module
is Gorenstein flat if  and only if it is Gorenstein projective. In
\cite[Proposition 1.3]{BM}, it is remarked that \cite[Theorem
5.1.11]{LW} can be generalized to coherent rings and for finitely
presented modules. Next we generalize this by showing, over
associative rings, that the same equivalence holds for infinitely
presented modules.\\
Recall that a left $R$-module $M$ is said to be infinitely
presented, if it admits a free resolution  $$\cdots  \rightarrow
F_{ 1} \rightarrow F_0\rightarrow M \rightarrow 0$$ such that each
$F_i$ is a finitely generated free left $R$-module. For instance,
over a left Noetherian ring, every finitely generated left module
is infinitely presented; and generally,  over a left coherent
ring, every finitely presented left module is infinitely
presented.

\begin{thm} \label{thm-Gp-Gf} Let $M$ be an infinitely presented left
$R$-module. Then,  $M$ is Gorenstein flat if and only if  it is
Gorenstein projective.
\end{thm}

The proof of this theorem involves the following lemma which
generalizes \cite[Lemma 5.1.10]{LW}.\bigskip

Recall that an exact sequence of finitely generated free left
$R$-modules $\mathbf{L}$ is called a complete resolution by
finitely generated  free left $R$-modules, if the dual complex $
\Hom_R (\mathbf{L}, R) $ is exact \cite[Definition 4.1.2]{LW}.

\begin{lem} \label{lem-Gp-Gf} Let $\mathbf{L}$ be  an exact sequence of finitely generated free
left $R$-modules. The following are equivalent:\begin{enumerate}
    \item $\mathbf{L}$ is a complete resolution by  finitely generated  free
    left $R$-modules;
    \item  $\mathbf{L}$ is a complete projective resolution;
    \item $\mathbf{L}$ is a complete flat resolution.
\end{enumerate}
\end{lem}
\proof First, using \cite[Lemma 3.59]{Rot}, the equivalence $(1)
\Leftrightarrow(3)$ has the same proof as the that of  $(i)
\Leftrightarrow(iii)$ of \cite[Lemma 5.1.10]{LW} (see errata on the
Christensen's homepage:
http://www.math.ttu.edu/$\sim$lchriste/).\\
\indent Now, since (2) is stronger than (1), it remains to prove the
implication $(1) \Rightarrow(2)$.\\ Let $$ \mathbf{L}=\quad
\cdots\rightarrow L_1\rightarrow L_0 \rightarrow L_{-1} \rightarrow
L_{-2} \rightarrow\cdots$$ be a complete resolution by finitely
generated free left $R$-modules. We decompose $\mathbf{L}$ into
short exact sequences
$$(\alpha_i)=\qquad 0\rightarrow N_{i +1}\rightarrow L_i\rightarrow N_{i} \rightarrow
0,$$ where  $N_i=\Im(L_i\rightarrow L_{i-1})$ for $i\in \Z$. To
prove that $\mathbf{L}$  is a complete projective resolution, it
is sufficient to prove, for every $i\in \Z$, that
$\Ext^{1}_R(N_i,Q)=0$ for every projective left $R$-module $Q$.
Indeed, if such condition holds, then the sequences
$$0\rightarrow \Hom(N_{i},Q)\rightarrow \Hom(L_{i},Q)\rightarrow
\Hom(N_{i +1},Q) \rightarrow 0$$ are all exact for every
projective left $R$-module $Q$. This implies, by assembling the
sequences $(\alpha_i)$, that $\mathbf{L}$   is a complete
projective resolution.\\
Then, we prove, for every $i\in \Z$, that $\Ext^{1}_R(N_i,Q)=0$
for every projective left $R$-module $Q$. Note first that each
$R$-module $N_i$ is infinitely presented and satisfies
$\Ext^{n}_R(N_i,R)=0$ for all $n>0$.\\ By Lazard's Theorem
\cite[\S 1, N$^\mathrm{o}$ 6, Theorem 1]{Bou}, there exists, for
every flat left $R$-module $F$, a direct system $(F_{j})_{j\in J}$
of finitely generated free left $R$-modules over a directed index
set $J$ such that $ \underrightarrow{lim}\;F_{j}\, \cong \, F$.
From \cite[Exercise 3, p. 187]{Bou} (or similarly to the proof of
\cite[Lemma 3.1.16]{Rel-hom}), we get, for every $i\in \Z$,
\begin{eqnarray*}
\Ext^{1}_R(N_i,F)&\cong&\Ext^{1}_R(N_i,\underrightarrow{lim}\;F_{j}) \\
         &\cong&\underrightarrow{lim}\;\Ext^{1}_R(N_i,F_{j}).
\end{eqnarray*}
Now, since every direct sum is the direct limit of its finite
partial sum ordered by inclusion, \cite[Exercise 3, p. 187]{Bou}
implies also that $\Ext^{1}_R(N_i,F_j)=0$. Therefore, for every
projective (then flat) left $R$-module $Q$, $\Ext^{1}_R(N_i,Q)=0$,
which completes the proof. \hfill$\square$\bigskip

\noindent\textit{Proof of Theorem \ref{thm-Gp-Gf}.} Using Lemma
\ref{lem-Gp-Gf} above, the ``if'' part is proved along the same
lines as \cite[Theorem 4.2.6]{LW} and the ``only if'' part is
proved similarly to \cite[Theorem 5.1.11]{LW}.
\hfill$\square$\bigskip

There are examples, over Noetherian rings, of finitely generated
(then infinitely presented) modules which are Gorenstein
projective (then   Gorenstein flat) but they are not projective
(then not flat) (see for instance \cite[Examples 1.1.13 and
4.1.5]{LW}). Namely, any ideal of a 1-Gorenstein ring (i.e.,
Noetherian with self-injective dimension at most 1) is infinitely
presented and Gorenstein projective. \bigskip

We take advantage of Lemma \ref{lem-Gp-Gf} and its proof to give a
result for the Gorenstein projective dimension. Precisely, using
the proof of Lemma \ref{lem-Gp-Gf}, we easily deduce the following
characterization of Gorenstein projective dimension of an
infinitely presented left module. This generalizes \cite[Theorem
1.2.7]{LW} (see also \cite[Theorem 4.4.12 and Corollary
4.4.13]{LW}).

\begin{cor}\label{cor-infinty-Gpd}
 Let $M$ be an infinitely presented left $R$-module with
finite Gorenstein projective dimension and let $n\geq 0$ be a
positive integer. Then, the following are equivalent:
\begin{enumerate}
    \item $\Gpd_R(M) \leq n$;
    \item $\Ext^{i}_R(M,F)=0 $ for all  $i>n$ and all flat left
    $R$-modules $F$;
    \item $\Ext^{i}_R(M,F)=0 $ for all  $i>n$ and all left
    $R$-modules $F$ of finite flat dimension;
    \item $\Ext^{i}_R(M,R)=0 $ for all  $i>n$.
\end{enumerate}
Consequently, the Gorenstein projective dimension of $M$ is also
determined by the formulas:
\begin{eqnarray*}
  \Gpd_R(M) &=& \sup\{ i\in \N \,| \, \Ext^{i}_R(M,F)\not = 0 \ for \ some\ flat\ left\ R\!-\!module \ F\} \\
            &=& \sup\{ i\in \N \,| \, \Ext^{i}_R(M,F)\not = 0 \ for \ some\ left\ R\!-\!module \ F\ with\ \fd_R(F)<\infty\} \\
  &=& \sup\{ i\in \N \,|\, \Ext^{i}_R(M,R)\not = 0
  \}.
\end{eqnarray*}
\end{cor}

It is important to note that the condition ``$M$ has finite
Gorenstein projective dimension'' in Corollary
\ref{cor-infinty-Gpd} above can not be dropped. Indeed, Jorgensen
and \c{S}ega \cite[Theorem 1.7]{JS} constructed, over an artinian
ring $R$, a finitely generated (then infinitely presented)
$R$-module  $M$  such that $\Gpd_R(M)=\infty$ and
$\Ext^{i}_R(M,R)=0 $ for all $i>0$.\bigskip

\noindent \textit{Acknowledgements}. The author thanks the referee
for very  helpful comments and suggestions.

\end{section}



\begin{thebibliography}{999}\addcontentsline{toc}{section}{\protect\numberline{}{Bibliography}}

\bibitem{BM} D. Bennis and N. Mahdou, Strongly Gorenstein projective, injective, and flat modules, \textit{J. Pure Appl. Algebra} \textbf{210} (2007) 437--445.

\bibitem{Bou} N. Bourbaki, \textit{Algèbre Homologique}, Chapitre 10, Masson, Paris, 1980.

\bibitem{CD} J. Chen and N. Ding, The flat dimensions of injective modules, \textit{Manuscripta Math.} \textbf{78} (1993) 165--177.

\bibitem{FCring} J. Chen and N. Ding, Coherent rings with finite self-FP-injective dimension, \textit{Comm. Algebra} \textbf{24} (1996) 2963--2980.

\bibitem{LW} L. W. Christensen, \textit{Gorenstein dimensions}, Lecture Notes in Math., 1747, Springer, Berlin,  2000.

\bibitem{CFH} L. W. Christensen, A. Frankild, and  H. Holm, On Gorenstein projective, injective and flat dimensions - a functorial description with applications,  \textit{J. Algebra} \textbf{302} (2006) 231--279.

\bibitem{copure} E. E. Enochs and O. M. G. Jenda,  Copure injective resolutions, flat resolvents and dimensions, \textit{Comment. Math. Univ. Carolin.} \textbf{34}   (1993)  203--211.

\bibitem{Rel-hom} E. E. Enochs and O. M. G. Jenda,  \textit{Relative homological algebra}, de Gruyter Expositions in Mathematics, vol. 30, Walter de Gruyter \& Co., Berlin, 2000.

\bibitem{GoPlat} E. E. Enochs, O. M. G. Jenda, and   B. Torrecillas,  Gorenstein flat modules, \textit{Nanjing Daxue Xuebao Shuxue Bannian Kan} \textbf{10} (1993) 1--9.

\bibitem{Fox}   E. E. Enochs, O. M. G. Jenda, and J. Xu, Foxby duality and Gorenstein injective and projective modules, \textit{Trans. Amer. Math. Soc.} \textbf{348}  (1996) 3223--3234.

\bibitem{Gf-cover} E. Enochs and J. Xu, Gorenstein Flat Covers of Modules over Gorenstein Rings, \textit{J. Algebra} \textbf{181} (1996) 288--313.

\bibitem{HH} H. Holm, Gorenstein homological dimensions, \textit{J. Pure Appl. Algebra} \textbf{189} (2004) 167--193.

\bibitem{JS} David A. Jorgensen and Liana M. \c{S}ega, Independence of the total reflexivity conditions for modules, \textit{Algebr. Represent. Theory} \textbf{9} (2006) 217--226.

\bibitem{Rot} J. Rotman,  \textit{An Introduction to Homological Algebra}, Pure and Applied Mathematics, 85, Academic Press, New York-London, 1979.

\end{thebibliography}
\end{document}